\documentclass[11pt]{article} 
\usepackage[english]{babel}
 \usepackage{amssymb,amsbsy,amsmath,amsfonts}

\usepackage{graphics,graphicx} 
\usepackage{color} 
\usepackage{graphicx}
\def\R{{\mathbb R}} 
  
\def\v3{\vskip0.3cm \noindent}

 \newtheorem{rem}{\bf Remark}

\newtheorem{thm}{\bf Theorem\/}

\begin{document} 
\date{today} 
\centerline{\Large {\bf Blow up of the solutions to a linear elliptic  }} 
\centerline{\Large {\bf system  involving Schr\" odinger operators}} 
\vskip0.3cm \noindent 
\centerline{\bf B.Alziary}
\centerline{\bf J.Fleckinger}
\centerline{\bf Institut de Math\'ematique, UMR5219}
\centerline{\bf Ceremath - Universit\'e Toulouse 1}
\vskip0.3cm \noindent
{\bf Abstract} $\;$ 
We show how   the solutions to a $2\times 2$  linear system involving Schr\" odinger operators 
blow up as the parameter $\mu$ tends to some critical value which is the principal eigenvalue of the system; 
here  the potential is  continuous  positive with superquadratic growth   and the square matrix of the system is  with constant coefficients   and may  have a double eigenvalue.
\vskip0.5cm \noindent
\section{Introduction} 
We study here the behavior  of the solutions to a $2\times 2$   system  (considered in its variational formulation):
$$LU:=(-\Delta + q(x))U = AU+\mu U + F(x)\; in \; \R^N,  \leqno(S)$$
$$U(x)_{|x|\rightarrow \infty} \rightarrow 0 $$
where $q$ is a continuous  positive potential tending to $+\infty$ at infinity  with superquadratic growth; $U$ is a column vector with components $u_1$ and $u_2$ and $A$ is a $2 \times 2$ 
square matrix with constant  coefficients.
$F$ is a column vector with components $f_1$ and $f_2$.
\par \noindent
Such systems have been intensively studied mainly for $\mu = 0$ and for $A$ with 2 distinct eigenvalues; here we consider also the case of a double eigenvalue. In both cases,  we show the blow up of solutions as $\mu$  tends to some critical value $\nu$ which is the principal eigenvalue of System  $(S)$.
This extends to systems involving Schr\" odinger operators defined on $\R^N$ earlier results valid for systems involving the classical Laplacian defined on smooth bounded domains with Dirichlet boundary conditions. 
\par \noindent
This paper is organized as follows: In Section 2 we recall known results for one equation. 
In Section 3 we consider first  the case where $A$ has two different eigenvalues and then  we study the case of a double eigenvalue. 

\par \noindent


\section{The equation}
We shortly recall the case of one equation 
$$ Lu:= (-\Delta +q(x)) u \, = \, \sigma u + f(x) \in \R^N, \leqno(E)$$ 
$$ \lim_{|x| \Rightarrow + \infty}  u(x) \, = \, 0.$$
$\sigma$ is a real parameter. 
\vskip0.3cm \noindent
{\bf Hypotheses}
\par \noindent
$(H_q)$ $q$ is a positive continuous potential tending to $+\infty$ at infinity.
\par \noindent
$(H_f)$ $\,$ $f \in L^2(\R^N)$, $f \geq 0$ and $f>0$ on some  subset with positive Lebesgue measure.
\par \noindent 
It is well knwon that if $(H_q)$ is satisfied,  $L$ possesses an infinity of 
eigenvalues tending to $+\infty$: $ \, 0 < \lambda_1 < \lambda_2 \leq \ldots .$
\vskip0.3cm \noindent
{\bf Notation: $(\Lambda, \phi)$} $\,$ Denote by 
 $\Lambda$ the smallest eigenvalue of $L$;  it  is positive and simple  and  denote by $\phi$ 
the associated eigenfunction, positive and with $L^2$-norm $\|\phi\|=1$. 
\vskip0.3cm \noindent
It is classical   (\cite{EdEv}, \cite{ReSi})  that if $f>0$ and $\sigma <\Lambda$ the positivity is improved,
 or in other words, the maximum principle {\bf  (MP)} is satisfied:
$$f\geq 0,\, 
 \not \equiv 0\; \Rightarrow\;  u>0. \leqno(MP)$$ 
Lately, for potentials growing fast enough (faster than the harmonic oscillator),  another notion has been
 introduced (\cite{AFTa1999}, \cite{AFTa2001}, \cite{ATa2007}, \cite{ATa2009}) which improves the maximum (or antimaximum principle):  the "groundstate positivity" ({\bf GSP}) 
(resp. " negativity"  ({\bf GSN})) which means that there exists $k>0$ such that
\par \noindent
\centerline{$u> k \phi$ (GSP)  (resp. $u<-k \phi$ (GSN))}.  
\par \noindent
We also say shortly  "fundamenal positivity" or" negativity", or also "$\phi$-positivity" or "negativity". 
\par \noindent
The first steps in this direction use  a radial potential. Here we consider a small perturbation of a radial one as in 
  \cite{ATa2007}.
\vskip0.3cm \noindent
 {\bf The potential $q$} $\,$ We define first 
 a class   ${\cal P}$  of radial potentials:
\begin{equation} \label{P}
{\cal P}:=\{Q \in {\cal C}(\R_+,(0,\infty)) /\,  \exists R_0>0, Q'>0 \, a.e. \, on \, [R_0, \infty), \,  
\int_{R_0}^{\infty} Q(r)^{-1/2} < \infty\}.\end{equation}
The last inequality  holds if $Q$ is growing sufficiently fast ($>r^2$).  Now we give  results of GSP or GSN for a 
 potential $q$ which is  a small perturbation of $Q$; we assume: 
\par \noindent
$(H'_q)$ $\,$ $q$ satisfies $(H_q)$ and there exists two functions $Q_1$ and $Q_2$ in ${\cal P}$, and 
two positive constants $R_0$ and $C_0$ such that 
\begin{equation} \label{qQ1}
Q_1(|x|) \leq q(x) \leq Q_2(|x|)  \leq C_0 Q_1 (|x|), \; \forall x\in \R^N, \end{equation} 
\begin{equation} \label{qQ2}
\int_{R_0}^{\infty} (Q_2(s) - Q_1(s)) \int_{R_0}^s 
exp\big( - \int_r^s [Q_1(t)^{1/2} + Q_2(t)^{1/2}]dt \big) drds <\infty.
\end{equation} 
Denoting by $\Phi_1$ (resp. $\Phi_2$) the groundstate of $L_1:= - \Delta + Q_1$ (resp. $L_2= - \Delta + Q_2$),  Corollary 3.3 in \cite{ATa2007} says that all these groundstates are "comparable" that is there exists constants $0<k_1 \leq k_2 \leq \infty$ such that $k_1\phi \leq \Phi_1, \Phi_2 \leq k_2\phi$. 
Finally 
\begin{thm}\label{gsp} (GSP)  (\cite{ATa2007}) $\,$
If  $(H'_q)$ and $(H_f)$ are satisfied, then, for $\sigma < \Lambda$,
there is a unique solution $u$ to $(E)$ which is positive, and there exists a constant $c>0$, such that 
\begin{equation} \label{min}     u>c\phi.\end{equation}
Moreover, if also $f\leq C\phi$ with some constant $C>0$, then 
\begin{equation} \label{maj}
u \leq \frac{C}{\Lambda - \sigma} \phi.\end{equation}
\end{thm}
{\bf The space ${\cal X}\, $} It is convenient for several results  to introduce the space of "groundstate bounded functions": 

\begin{equation} \label{X}
{\cal X}:=\{ h \in L^2(\R^N): \, h/\phi \in L^{\infty}(\R^N)\}, \end{equation}
equipped with the norm $\|h\|_{\cal X}=ess\sup_{\R^n}(|h|/\phi)$.
\par \noindent
For  a potential satisfying $(H'_q)$ and a function $f \in {\cal X}$, there is also a result of "groundstate negativity" {\bf (GSN)} for $(E)$; it is  is an extension of the antimaximum principle, introduced by Cl\'ement and Peletier in 1978 (\cite{ClPe}) 
for the Laplacian when the parameter $\sigma$ crosses $\Lambda$. 
\begin{thm}\label{GSN} (GSN) (\cite{ATa2007} ) $\,$  Assume  $(H'_q)$ and $(H_f)$ are satisfied and $f\in {\cal X}$; then there exists $\delta(f)>0$ and a positive constant $c'>0$ such that for all
$\sigma \in (\Lambda, \Lambda + \delta)$, 
\begin{equation} \label{neg} u \leq -c'\phi. \end{equation}
\end{thm}
\begin{rem}
This holds also if we only assume $f^1:= \int f\phi>0.$
\end{rem}
\par \noindent
{\bf
Hypothesis $(H'_f)$}$\,$ We consider  now  functions $f$ which are such that
\par \noindent
$(H'_f)$: $\,$ $f \in {\cal X}$ and $f^1:= \int f \phi >0$.
\par \noindent
\begin{thm}\label{LIM} $\,$ Assume $(H'_q)$ and $(H'_f)$ are satisfied. Then there exists $\delta >0$ such that for 
 $\Lambda - \delta <  \sigma < \Lambda$ there exists  positive constants $k'$ and  $K'$,  depending on $f$ and $\delta$ such that
\begin{equation} \label{MIN}
0 < \frac {k'}{\Lambda - \sigma } \phi < u <  \frac {K'}{\Lambda - \sigma } \phi.
\end{equation}
If   $ \Lambda  < \sigma < \Lambda + \delta$,  there exists  positive constants $k"$ and $K"$,  depending on $f$ and $\delta$ such that
\begin{equation} \label{MAX}
 - \frac {k"}{\Lambda - \sigma } \phi <   - u < - \frac {K"}{\Lambda - \sigma } \phi  <0. 
\end{equation}
\end{thm}
This result extends earlier one in \cite{Le2010}  and a a close  result is Theorem 2.03  in \cite{Be}. It shows in particular that $u \in {\cal X}$ and $|u| \rightarrow \infty$ as $|\nu - \mu| \rightarrow 0$. 
\par \noindent
{\bf Proof:} $\,$ 
Decompose  $u$ and $f$  on $\phi$ and its orthogonal:
\begin{equation}\label{DECOMP} u= u^1 \phi +  u^{\perp}\, ; \; f=f^1 \phi + f^{\perp}.\end{equation}
We derive from $(E)$  $Lu = \sigma u + f$ : 
\begin{equation}\label{PERP} L u^{\perp} = \sigma u^{\perp} + f^{\perp}\, \end{equation}
\begin{equation}\label{E1}  L u^1 \phi = \Lambda u^1 \phi = \sigma   u^1 \phi + f^1 \phi.\end{equation}
We notice that since $q$ is smooth; so is $u$. Also, since $f \in {\cal X}$, $f^{\perp}$,  $u$ and  $u^{\perp}$ are also in 
 ${\cal X}$ and hence are bounded. 
Choose $ \sigma < \Lambda$ and assume $(H'_f)$. 
We derive from Equation (\ref{PERP})  (by \cite{AFTa2007}Thm 3.2)  that : $||u^{\perp}||_{ {\cal X}} < K_1$.   Therefore $|u^{\perp}| $ is bounded by some $cste. \phi>0$.
\par \noindent
From Equation  (\ref{E1}) we derive
 \begin{equation}\label{u1} u^1= \frac{f^1}{(\Lambda - \sigma)} \, \rightarrow \, \pm \infty \, as
\, (\Lambda - \sigma)  \rightarrow 0.\end{equation}
Choose $\Lambda - \delta $  small enough and $\sigma \in (\Lambda - \delta, \Lambda)$.
 Hence
$$
\frac {K'}{\Lambda - \sigma } \phi < u;  \;\; u <  \frac {K"}{\Lambda - \sigma } \phi.
$$
For $\sigma > \Lambda$. we do exactly the same, except that the signs are changed for $u^1$ in $(\ref{u1})$.

\vskip0.3cm \noindent
\section{ A $2\times 2$  Linear  system}
Consider now a linear system with constant coefficients.
$$ LU=AU+\mu U + F(x) \; in \; \R^N. \leqno(S)$$
As above, $L:= - \Delta + q$ where the potential $q$ satisfies $(H'_q)$, and where $\mu$ is a real parameter. 
\noindent $L$ can be detailed as 2 equations:
$$ \left \{ \begin{array}{lcl}
Lu_1&=&au_1+bu_2+ \mu u_1 +f_1(x)\\
 Lu_2&=&c u_1+du_2+ \mu u_2 +f_2(x)\end{array} \right .\, \; in \; \R^N ,.\leqno(S)$$
$$ \; u_1(x), u_2(x)_{|x|\rightarrow \infty} \rightarrow 0 .$$
Assume 
$$A=\left(\begin{array}{cccc}
a&b\\c&d
 \end{array}
 \right)\,  \mbox{ with } \, b>0\,  {\rm and}\, D:= (a-d)^2+4bc \geq 0. \leqno(H_A)$$
Note that $b>0$ does not play any role since we can always change the order of the equations. 
\par \noindent
The eigenvalues of $A$  are $$\xi_1= \frac{a+d + \sqrt D}{2} \geq \xi_2= \frac{a+d - \sqrt D}{2} .$$
As far as we know, all the previous studies suppose that  the largest eigenvalue  $\xi_1$ is simple (i.e. $D=(a-d)^2+4bc>0$).
Here we  also study, in the  second subsection,  the case  of  a double eigenvalue $\xi_1 = \xi_2$, that is  $D=0$; this implies necessarily $bc<0$ and necessarily the matrix is not cooperative.
\subsection{Case $\xi_1>\xi_2$}$\,$ This is the classical case where $\xi_1$ is simple. 
Set $\xi_1>\xi_2$.
The eigenvectors are 
$$X_k= \left(\begin{array}{cccc}
b \\  \xi_ k  -a 
 \end{array}
 \right),$$
Set  $X:=X_1$. 
\par \noindent
As above, denote by  $(\Lambda , \phi)$,\quad  $\phi>0$,  the  principal eigenpair
of  the operator  $L=(-\Delta + q(x))$.  
\par \noindent
It is easy to see that 
$$ L(X\phi)  - A X \phi = (\Lambda -   \xi_1 )  X\phi.$$
Hence \begin{equation}\label{NU}\nu=\Lambda - \xi_1\end{equation}  is the principal eigenvalue of $(S)$ with associated eigenvector $X\phi$.
Note that the components of $X\phi$ do not change sign, but, in the case of a non cooperative matrix 
they are not necessarily both positive.
We prove:

\begin{thm}\label{THSS} $\,$ Assume $(H'_q)$and $(H_A)$;   $f_1$ and $f_2$ satisfy $(H'_f)$; assume also $D>0$ and  $d-a>0$. If 
$$(\xi_2 - a )f_1^1< bf_2^1, $$ there exists $\delta >0$, independant of $\mu$, 
such that if $\nu- \delta < \mu < \nu$, there exists a positive constant $\gamma$  depending only on $F$ such that 
\begin{equation} \label{THS1}
u_1, u_2 \geq \frac{\gamma}{\nu - \mu} \phi >0.\end{equation} 
\par \noindent
If  $\nu < \mu < \nu+  \delta $, the sign are reversed:
\begin{equation} \label{THS2}
u_1, u_2 \leq - \frac{\gamma}{\nu - \mu} \phi <0.\end{equation} 
\end{thm}
\begin{rem} If $(H'_q)$ and $(H_A)$ are satisfied; if   $f_1$ and $f_2$ satisfy $(H'_f)$ as in Theorem \ref{THSS}, but 
 if  $d-a<0$ we have, if $(a-\xi_2) f_1^1 + bf_2^1>0$  and $\nu- \delta < \mu < \nu$
$$u_1 \geq \frac{\gamma}{\nu - \mu} \phi >0, \;  u_2 \leq  - \frac{\gamma}{\nu - \mu} \phi <0 $$
\end{rem} 
\begin{rem} It is noticeable that for all these cases, $|u_1|$, $|u_2|$ $\rightarrow +\infty$ as $|\nu - \mu| \rightarrow 0$.
\end{rem} 
These results  extend Theorem 4.2 in \cite{AFTa1999}.

\par \noindent
{\bf Proof: }$\,$ 
As in \cite{AF2016}, we use $J$ the associated Jordan matrix (which in this case is diagonal) 
 and $P$ the change of basis matrix which are  such that 
$$A=PJP^{-1}.$$
Here
\begin{equation} \label{matP} P =\left(\begin{array}{cccc}
b&b\\\xi_1 - a &\xi_2 - a
 \end{array}
 \right),\quad P^{-1} = \frac{1}{b(\xi_1 - \xi_2)}\left(\begin{array}{cccc}
a-\xi_2&b\\\xi_1 - a & -b
 \end{array}
 \right).\end{equation}
$$J=\left(\begin{array}{cccc}
\xi_1&0\\0&\xi_2
 \end{array}
 \right).$$
Denoting $\tilde U = P^{-1} U$ and  $\tilde F = P^{-1} F$, we derive from System  $(S)$ (after multiplication by  $\ P^{-1} U$
to the left):
$$ L \tilde U = J \tilde U + \mu \tilde U + \tilde F.$$
Since $J$ is diagonal we have two independant equations:

\begin{equation}\label{DIAG} L \tilde u_k = ( \xi_k + \mu)  \tilde u_k + \tilde f_k, \; k=1\, or \, 2.\end{equation}
The projection on $\phi$ and on its orthogonal for $k=1$ and $2$ gives 
$${ \tilde u}_k =( { \tilde u}_k)^1 \,\phi  +  {  \tilde u}_k^{\perp}, \quad { \tilde f}_k = ({ \tilde f}_k)^1 \,\phi  +
 { \tilde f}_k^{\perp};$$
hence
\begin{equation}\label{tildeEk}
L (\tilde u_k)^1 \,\phi = \Lambda ( \tilde u_k)^1 \,\phi =  \xi_k (\tilde u_k)^1 \,\phi + \mu (\tilde u_k)^1 \,\phi +( \tilde f_k)^1 \phi,\end{equation}
\begin{equation}\label{tildeE'k}
L {\tilde u_k}^{\perp} = \xi_k {\tilde u}_k^{\perp} + \mu {{\tilde u}_k}^{\perp} + {\tilde f_k}^{\perp}.\end{equation}
If both $f_k$ verify $(H'_f)$ , they are are in ${\cal X}$ and bounded and hence both ${ \tilde f}_k^{\perp}$ are bounded; therefore, by  (\ref{tildeE'k})
both $ { \tilde u}_k^{\perp}$ are also bounded. 
\par \noindent
We derive from (\ref{tildeEk}) that 
$$(\tilde u_k)^1\,  =\,  \frac{(\tilde f_k)^1}{\Lambda - \xi_k - \mu} \, = \, \frac{(\tilde f_k)^1}{\nu - \mu}.$$ 
Consider again  Equation  (\ref{tildeEk}) for $k=2$;  obviously, $ (\tilde u_2)^1$ stays bounded as $\mu \rightarrow \nu = \Lambda -  \xi_1$ and therefore $\tilde u_2$ stays bounded. . 
\par \noindent
For $k=1$, $(\tilde u_1)^1 = \frac{(\tilde f_1)^1}{\nu - \mu} \rightarrow \infty$ as 
$\mu \rightarrow \nu = \Lambda -  \xi_1$, since $(\tilde f_1)^1=\frac{1}{\xi_1 - \xi_2}((a - \xi_2  )f_1^1+ bf_2^1 )>0$,; this is the condition which appears in Theorem \ref{THSS}. 
Then, we simply apply Theorem \ref{LIM} to (\ref{DIAG}) for $k=1$  and deduce that 
 there existes $\delta >0$,  such that, for $|\Lambda -\xi_1-\mu|=|\nu - \mu|< \delta$, there exists a positive constant $C>0$ such that $\tilde u_1 \phi \geq \frac{C}{\nu -  \mu}\phi$.
Now, it follows from
 $U=P \tilde U$, that $$u_1=b(\tilde  u_1 + \tilde u_2), \; u_2 =  (\xi_1 - a) \tilde u_1 +  (\xi_2 - a) \tilde u_2.$$
As $\nu - \mu \rightarrow 0$, since $\tilde u_2$ stays bounded,  $u_1$ behaves as  $ b(\tilde u_1)^1\phi$, 
 $ u_2 $ as $ (\xi_1 - a ) (\tilde u_1)^1\phi $. 
More precisely, if $|\mu - \nu|$ small enough
$$(\tilde u_1)^1 \geq  \frac{K}{\nu-\mu}   \; if \; \mu < \nu \,; \; \tilde u^1_1 \leq  - \frac{K}{\nu-\mu}   \; if \; \mu > \nu  $$
where $K$ is a positive constant depending only on $F$.

\begin{rem}
Indeed, we always assume that $b>0$, hence $u_1>0$ for $\nu - \mu>0$ small enough.  
For the sign of $u_2$ we remark  that  $(\xi_1 - a)$ and $(d-a)$ have the same sign. 
\end{rem}

\subsection{Case $\xi_1=\xi_2$}
Consider now the case where the coefficients of the matrix $A$ satisfy  $(H_A)$ and 
$$D:= (a-d)^2+4bc=0.$$
\par \noindent
Of course this implies $bc<0$ and since  $b>0$ , then  $c<0$:  we have a non cooperative system. Now $\xi_1=\xi_2=\xi=\frac{a+d}{2}$.
We prove here
\begin{thm}\label{THSD}Assume $(H'_q)$ and $(H_A)$ with  $(a-d)^2+4bc=0$; assume also that  $f_1 ,f_2$   satisfy $(H'_f)$ and :
$$ \frac{(a-d)}{2}f_1^1+ bf_2^1>0.$$ 
  If $ \mu < \nu =  \Lambda - \xi$, $\nu - \mu < \delta$, small enough,  there exists a positive constant $\gamma$ such that 
$$u_1\geq\frac{\gamma}{\nu - \mu} \phi, \;\;  
u_2 \leq -\frac{\gamma}{\nu - \mu} \phi.$$ 
 If $  \nu =  \Lambda - \xi < \mu < \nu + \delta$, ($\delta$ small enough),  there exists a positive constant $\gamma'$ such that 
$$u_1\leq - \frac{\gamma}{\nu - \mu} \phi, \; \; 
u_2 \geq \frac{\gamma}{\nu - \mu} \phi.$$ 
\end{thm} 
\par 
\begin{rem}
Note that the condition $ \frac{(a-d)}{2}f_1^1+ bf_2^1>0$ in the theorem above is the same than in theorem \ref{THSS} $(\xi_2 - a )f_1^1< bf_2^1$, since in theorem \ref{THSD} $\xi_2=\xi=\frac{a-d}{2}$.
\end{rem}
\noindent
{\bf Prrof} $\,$  The eigenvector  associated to eigenvalue $\xi$ is 
$$X= \left(\begin{array}{cccc}
b \\  \frac{d  -a }{2}
 \end{array}
 \right).$$
The vector $X\phi$ is thus an eigenvector for $L-A$,
$$L(X\phi)-AX\phi=(\Lambda -\xi)X\phi = \nu X \phi.$$
We use again $J$ the associated Jordan matrix 
and $P$ the change of basis matrix; we have 
$$A=PJP^{-1}.$$
Here
$$P =\left(\begin{array}{cccc}
b&\frac{2b}{a-d}\\\frac{d- a }{2}&0
 \end{array}
 \right),\quad P^{-1} = \frac{1}{b}\left(\begin{array}{cccc}
0&-\frac{2b}{a-d}\\\frac{a-d}{2}  & b
 \end{array}
 \right). $$
$$J=\left(\begin{array}{cccc}
\xi&1\\0&\xi
 \end{array}
 \right).$$
As above, setting $\tilde U = P^{-1} U$ and  $\tilde F = P^{-1} F$, we derive from System  $(S)$
$$ L \tilde U = J \tilde U + \mu \tilde U + \tilde F.$$
We do not have anymore a decoupled system but 
\begin{equation}\label{TILD}\left \{ \begin{array}{lclcl}
L\tilde{u}_1&=&(\xi +\mu)\tilde{u}_1&+&\;\, \tilde{u}_2+\tilde{f}_1\\
L\tilde{u}_2&=&&+&(\xi+\mu) \tilde{u}_2+\tilde{f}_2
\end{array}\right .\end{equation}
If  $\xi + \mu < \Lambda$ (that is $ \mu < \nu$) and if $\tilde{f}_2=\dfrac{(a-d)}{2b}f_1+f_2$ and $\tilde{f}_1=\frac{-2}{a-d}f_2$   satisfies $(H'_f)$, hence are in ${\cal X}$ and  $\frac{(a-d)}{2b}f_1^1+ f_2^1>0  $.  By Theorem \ref{LIM} applied to the second equation,   there exists a constant $K>0$, such that $\tilde{u}_2> \frac{K}{\nu - \mu} \phi$. Hence,  for $\nu-\mu$ small enough  fo any $ \tilde f_1 \in {\cal X}$, $\tilde{u}_2+\tilde{f}_1>0$ and is in $X$; then   again Theorem \ref{LIM} for the first equation implies that  there exists a constant $K'>0$, such that $\tilde{u}_1>\frac{K'}{\nu - \mu}\phi$.
\par \noindent
Since here  $a>d$., 
 there exists a  constant $\gamma >0$,
$$U=P\tilde{U}=\left \{ \begin{array}{l}
u_1 \, = \, b\tilde{u}_1+\frac{2b}{a-d}\tilde{u}_2\, > \,  \frac{\gamma}{\nu - \mu} \phi \\
\, u_2\, = \, \frac{d-a}{2}\tilde{u}_1\, < \, -  \frac{\gamma}{\nu - \mu}  \phi
\end{array}\right.$$

Again $u_1 \rightarrow +\infty$ as $\nu - \mu \rightarrow 0$ and  $u_2$ $\rightarrow -\infty$ as $\nu - \mu \rightarrow 0$.
\par
\noindent If $\mu>\Lambda$ ( and $\mu-\nu>0$ small enough we have analogous calculation with signs reversed. 
\begin{rem}
The results in theorem \ref{THSD} coincide with those of theroem \ref{THSS} in the case $D=0$.
\end{rem}

\color{black}

{}

\end{document}